\documentclass[conference]{IEEEtran}
\usepackage[utf8]{inputenc}
\usepackage[T1]{fontenc}

\usepackage{pgfplots}
\pgfplotsset{
    compat=newest,
    colormap={parula}{
        rgb255=(53,42,135)
        rgb255=(15,92,221)
        rgb255=(18,125,216)
        rgb255=(7,156,207)
        rgb255=(21,177,180)
        rgb255=(89,189,140)
        rgb255=(165,190,107)
        rgb255=(225,185,82)
        rgb255=(252,206,46)
        rgb255=(249,251,14)
    },
}

\usepackage{caption}
\usepackage{subcaption}

\usepackage{amsmath, amssymb, mathtools}
\usepackage{bm} 
\usepackage{dsfont} 

\newtheorem{theorem}{Theorem}[section]
\newtheorem{corollary}[theorem]{Corollary}

\newcommand{\field}[1]{\ensuremath{\mathds{#1}}}

\newcommand{\N}{\field N}

\newcommand{\C}{\field{C}}

\newcommand{\Z}{\field Z}

\newcommand{\R}{\field R}
\newcommand{\T}{\field T}

\newcommand{\bk}{\mathbf{k}}

\newcommand{\bu}{\mathbf{u}}
\newcommand{\bw}{\mathbf{w}}
\newcommand{\bx}{\mathbf{x}}

\newcommand{\be}{\begin{equation}}
\newcommand{\ee}{\end{equation}}
\newcommand{\beq}{\begin{eqnarray}}
\newcommand{\beqq}{\begin{eqnarray*}}
\newcommand{\eeq}{\end{eqnarray}}
\newcommand{\eeqq}{\end{eqnarray*}}

\newtheorem{remark}[theorem]{Remark}

\definecolor{azure}{RGB}{0, 127, 255}
\definecolor{violet}{RGB}{161, 11, 112}
\definecolor{orange}{RGB}{243, 119, 53}

\usepackage{todonotes}

\usepackage{orcidlink}

\title{Efficient recovery of non-periodic multivariate functions from few scattered samples}

\author{\IEEEauthorblockN{Felix Bartel\IEEEauthorrefmark{1},
Kai Lüttgen\IEEEauthorrefmark{2}, Nicolas Nagel\IEEEauthorrefmark{1} and
Tino Ullrich\IEEEauthorrefmark{1}}
\IEEEauthorblockA{Faculty of Mathematics,
Chemnitz University of Technology\\
09107 Chemnitz, Germany\\
Email: \IEEEauthorrefmark{1}\{felix.bartel, kai.luettgen, nicolas.nagel, tino.ullrich\}@mathematik.tu-chemnitz.de,
\IEEEauthorrefmark{2}luettgen.math@gmail.com}}

\date{\today}

\begin{document} 

\maketitle

\begin{abstract} It has been observed by several authors (see \cite{COOLS2016166, KMNN21, NP21c}) that well-known periodization strategies like tent or Chebyshev transforms lead to remarkable results for the recovery of multivariate functions from few samples. So far, theoretical guarantees are missing. The goal of this paper is twofold. On the one hand, we give such guarantees and briefly describe the difficulties of the involved proof. On the other hand, we combine these periodization strategies with recent novel constructive methods for the efficient subsampling of finite frames in $\field{C}^m$. As a result we are able to reconstruct non-periodic multivariate functions from very few samples. The used sampling nodes are the result of a two-step procedure. Firstly, a random draw with respect to the Chebyshev measure provides an initial node set. A further sparsification technique selects a significantly smaller subset of these nodes with equal approximation properties. This set of sampling nodes scales linearly in the dimension of the subspace on which we project and works universally for the whole class of functions. The method is based on principles developed by Batson, Spielman, and Srivastava \cite{BaSpSr09} and can be numerically implemented. Samples on these nodes are then used in a (plain) least-squares sampling recovery step on a suitable hyperbolic cross subspace of functions resulting in a near-optimal behavior of the sampling error. Numerical experiments indicate the applicability of our results.
\end{abstract}

\section{Introduction}

We investigate the worst-case approximation of non-periodic functions $f$ on the cube $[-1, 1]^d$ via function samples and propose a two-step procedure. Firstly, we have to find a good set of sample points and secondly, we need a method to use this information to get an approximation for functions from the target space.

In the periodic case, it is well-known that choosing sufficiently many uniform random sample nodes in $[-1, 1]^d$ together with a least squares recovery using the Fourier basis yields near-optimal reconstructions. The rate of convergence in this case is determined by the function model, in particular the smoothness $s$ (made rigorous in Section \ref{sec:spaces}) with polylogarithmic factors determining the precise behaviour in dimension $d$.

The non-periodic case on the other hand is not so clear. One may try to use a periodization technique to turn $f$ into a periodic function $\tilde{f}$ and apply the above described procedure to $\tilde{f}$. Extending for example $f$ periodically via $\tilde{f}(\bx) = f((\bx \mod 2) - \mathbf{1})$, i.e. ``copying'' $f$ over to other cubes, gives such a periodization but may also introduce discontinuities (as is the case for the function in Figure \ref{fig:Bspline2:nonper}). We thus aim to find a periodization, such that $\tilde{f}$ is at least as smooth as $f$ itself.

In Section \ref{sec:spaces} we give such a periodization via the cosine composition operator $\tilde{f} = T_{\cos} f$ acting between suitable Sobolev-type spaces by
\begin{align} \label{eq:CosComposition}
    (T_{\cos} f)(x_1, \dots, x_d) = f(\cos(\pi x_1), \dots, \cos(\pi x_d)).
\end{align}
For us the important consequence of this is that approximation algorithms for periodic functions can be taken over to the non-periodic case, keeping the cosine composition in mind.

If $U$ is a uniform random variable taking values in $[-1, 1]$, then $\cos(\pi U)$ is distributed according to the Chebyshev measure $(\pi\sqrt{1-x^2})^{-1} \mathrm{d}x$ on $[-1, 1]$. Similarly, since $\tilde{f} = T_{\cos} f$ is even, the $d$-dimensional Fourier basis (on $[-1, 1]^d$) boils down to $\prod_{\ell=1}^d \cos(\pi k_\ell x_\ell)$. Consequently, going over to the non-periodic case, we will approximate $f$ using $\prod_{\ell=1}^d \cos(k_\ell \arccos(x_\ell))$, i.e. tensor product Chebyshev polynomials.

The first part of the previous paragraph suggests that we need to sample points following the Chebyshev distribution
\begin{align} \label{eq:rho}
    \mathrm{d}\varrho^{}_D(\bx)=\prod_{\ell=1}^d\left(\pi\sqrt{1-x_\ell^2}\right)^{-1} \mathrm{d}\bx.
\end{align}
These random points are near-optimal.
To achieve the optimal error behaviour, we will use a subsampling technique to lower the number of used sample points even further while still keeping good approximation properties. This will be explained in Section \ref{sec:subsampling}.

Lastly, we will describe our method for function recovery in Section \ref{sec:numerics}. Our claims are backed up by numerical experiments with the test function $f$ from \eqref{eq:testfunction} (tensored cutouts of a quadratic B-spline), the two-dimensional version being depicted in Figure \ref{fig:Bspline2:nonper}. The resulting $L_2$-errors as shown in Figure \ref{fig:L2_error_cheb_subsampled} strongly follow a main decay rate of $s = 5/2$, the smoothness of $f$.

We have two aims with this paper. Firstly, we want to demonstrate the power and versatility of subsampling by applying it to the Chebyshev setting. And secondly, we want to give a theoretical explanation of why Chebyshev distributed sample points are suitable for the recovery of non-periodic functions. While other methods cap at certain decay rates or only give half the rate due to deterministic components (see Section \ref{sec:numerics} for such methods), our method gives an optimal main decay rate of $s$.


\section{Periodization Operator and Embeddings} \label{sec:spaces}

In preparation of what is to follow, we fix some notation and give a few definitions. $\N$ denotes the positive integers and $\N_{0} = \N \cup \{0\}$. The parameter $d \in \N$ is used for the dimension. Bold characters are used for vectors from some $d$-dimensional space. The real parameter $s$ is always assumed to satisfy $s > 0$. Let $D = [-1,1]^{d}$. On $D$ we distinguish between two $L_{2}$-spaces. Firstly, $L_{2}(D)$ denotes the usual space of square-integrable functions with respect to the Lebesgue measure. Secondly, $L_{2} \! \left(\varrho^{}_D\right)$ is the space of square-integrable functions on $D$ with respect to the Chebyshev measure from \eqref{eq:rho}. The inner product on this latter space is denoted by $\langle f, g \rangle_{\varrho}$.

The torus is $\T^d = (\R / (2\Z))^d$, where we will use $D$ as the fundamental domain. In other words, $\T^d$ is derived from $D$ by identifying opposite sides of $D$. $L_{2}(\T^{d})$ denotes the usual space of square-integrable periodic functions.

For $m \in \N_{0}$ the (classical) Sobolev space $H^{m}_{\textnormal{mix}}(D)$ of dominating mixed smoothness is the space of all $f \in L_{2}(D)$ such that $\lVert f \rVert_{H^{m}_{\textnormal{mix}}(D)} = \lVert f \rVert_{L_{2}(D)} + \sum_{\alpha} \lVert \partial^{\alpha} f \rVert_{L_{2}(D)} < \infty$, where we sum over all $\alpha \in \N_{0}^{d}$ with $0 \leq \max_{i \in [d]} \alpha_{i} \leq m$. In case of fractional smoothness $s > 0$, $H^{s}_{\textnormal{mix}}(D)$ is defined via the Besov-type norm
\begin{equation}\label{eq:besov_type_norm}
\textstyle
    \lVert f \rVert_{H^{s}_{\textnormal{mix}}(D)}
    = \lVert f \rVert_{L_{2}(D)} + \sum_{e \subset [d]} \mathcal{R}^{M}_{e}(f)
    < \infty.
\end{equation}
Here, $M \geq s$ is a fixed integer and $\mathcal{R}_{e}^{M}(f)$ is given by
\begin{equation}
\textstyle
    \mathcal{R}_{e}^{M}(f)^{2}
    = \sum_{\mathbf{j} \in \N_{0}^{d}(e)} 2^{2 s \lvert \mathbf{j} \rvert_{1}} \sup_{\lvert h_{\ell} \rvert \leq 2^{-j_{\ell}}} \lVert \Delta_{\mathbf{h}}^{m,e} f \rVert_{L_{2}(D)}^{2}, \notag
\end{equation}
where $\N_{0}^{d}(e) = \{ \mathbf{j} \in \N_{0}^{d} \; \vert \; j_{\ell} = 0 \; \forall \ell \notin e \}$, $\Delta_{\mathbf{h}}^{m,e} = \prod_{\ell \in e} \Delta_{h_{\ell}, \ell}^{m}$ and $\Delta_{h_{\ell},\ell}^{m}$ is the $m$th-order forward difference in the $\ell$th coordinate direction of step size $h_\ell$.

The periodic Sobolev space $H^{s}_{\textnormal{mix}}(\T^{d})$ consists of all $f \in L_{2}(\T^{d})$ for which
\begin{equation}
    \lVert f \rVert_{H^{s}_{\textnormal{mix}}(\T^{d})}^{2}
    = \sum_{\mathbf{k} \in \Z^{d}} \lvert \hat{f}_{\mathbf{k}} \rvert^{2} \prod_{\ell = 1}^{d} \left(1 + \lvert k_{\ell} \rvert^{2}\right)^{s}
    < \infty, \notag
\end{equation}
where the Fourier coefficients of a function $f$ are given by
$$
\hat{f}_{\mathbf{k}} = 2^{-d} \int_{\T^{d}} f(\mathbf{x}) \, \mathrm{e}^{-\pi \mathrm{i} \mathbf{k} \cdot \mathbf{x}} \, \mathrm{d} \mathbf{x}, \quad \mathbf{k} \in \Z^{d} .
$$

Finally, we consider the $L_{2}(\varrho_{D})$-normalized Chebyshev polynomials of the first kind
\begin{align} \label{eq:cheb_poly}
\textstyle
    T_{k}(x) = \sqrt{2}^{\min\{1, k\}} \cos(k \arccos(x)), \quad k \in \N_{0}.
\end{align}
Taking tensor products of these polynomials, we set
\begin{align} \label{eq:eta}
\textstyle
    \eta_{\mathbf{k}}(\mathbf{x}) = \prod_{\ell = 1}^{d} T_{k_{\ell}}(x_{\ell}), \quad \bk = (k_\ell)_{\ell=1}^d \in \N_0^d.
\end{align}

With these definitions in mind we move on to the topic of periodization. Of the many ways to extend functions in a periodic fashion we choose to study the operator $T_{\cos}$ as defined in \eqref{eq:CosComposition}. On the one hand, thinking purely in terms of smoothness, the analysis of $T_{\cos}$ comes down to the chain rule and its generalizations. On the other hand, $T_{\cos}$ is not bounded as an operator from $L_{2}(D)$ to $L_{2}(\T^{d})$, i.e. $T_{\cos}$ does not mesh well with integrability. Since Sobolev spaces measure a combination of both smoothness and integrability, the study of the action of $T_{\cos}$ on these spaces is somewhat delicate. Nonetheless, employing tools such as maximal functions, the Fefferman-Stein inequality and the Riesz-Thorin interpolation theorem we are able to prove the following theorem in \cite{LuUl23}.

\smallskip

\begin{theorem} \label{thm:Stetigkeit Kosinuskomposition}
If $s > 1/2$, the linear operator $T_{\cos}: H^{s}_{\textnormal{mix}}(D) \rightarrow H^{s}_{\textnormal{mix}}(\T^{d})$ is bounded.
\end{theorem}

\smallskip

There is another interesting property of the periodization operator $T_{\cos}$. A simple calculation reveals that $\left(T_{\cos}f\right)^{\land}_{\mathbf{k}} = \langle f, \eta_{\mathbf{k}} \rangle_{\varrho} \prod_{\ell = 1}^{d} \sqrt{2}^{-\min\{1, k_{\ell}\}}$. In light of the above theorem we conclude:

\smallskip

\begin{corollary}
If $s > 1/2$ and $f \in H^{s}_{\textnormal{mix}}(D)$, we have
\begin{equation}
    \sum_{\bk \in \N_{0}^{d}} \lvert \langle f, \eta_{\mathbf{k}} \rangle_{\varrho} \rvert^{2} \prod_{\ell = 1}^{d} \left(1 + \lvert k_{\ell} \rvert^{2}\right)^{s}
    \lesssim \lVert f \rVert_{H^{s}_{\textnormal{mix}}(D)}^{2}
    < \infty . \notag
\end{equation}
\end{corollary}

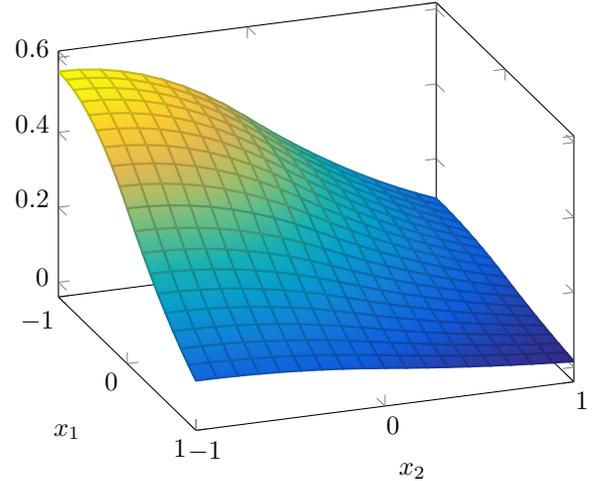
\begin{figure}[htb]
	\centering
	\begin{tikzpicture}
	
	\begin{axis}[
    view={70}{30}, 
    xlabel=$x_{1}$, ylabel=$x_{2}$,
	xtick={-1,0,1},
	ytick={-1,0,1},
    ]
    
    \addplot3[
    z buffer=sort,
	surf,
	domain=-1:0,
	domain y=-1:0,
	samples = 10,
    samples y = 10,
    shader = faceted interp,
    line width = .8pt,
    ] 
	{(-0.25*x*x - 0.5*x + 0.5)*(-0.25*y*y - 0.5*y + 0.5)};
	
	\addplot3[
	z buffer=sort,
	surf,
	domain=-1:0,
	domain y=0:1,
	samples = 10,
    samples y = 10,
    shader = faceted interp,
    line width = .8pt,
    ] 
	{(-0.25*x*x - 0.5*x + 0.5)*(0.125*y*y - 0.5*y + 0.5)};
	
	\addplot3[
	z buffer=sort,
	surf,
	domain=0:1,
	domain y=-1:0,
	samples = 10,
    samples y = 10,
    shader = faceted interp,
    line width = .8pt,
    ] 
	{(0.125*x*x - 0.5*x + 0.5)*(-0.25*y*y - 0.5*y + 0.5)};
	
	\addplot3[
	z buffer=auto,
	clip=false,
	surf,
	domain=0:1,
	domain y=0:1,
	samples = 10,
    samples y = 10,
    shader = faceted interp,
    line width = .8pt,
    ] 
	{(0.125*x*x - 0.5*x + 0.5)*(0.125*y*y - 0.5*y + 0.5)};
	
	\end{axis}
	
	\end{tikzpicture}
	\caption{Surface plot of the test function from (\ref{eq:testfunction}) for $d = 2$.}
	\label{fig:Bspline2:nonper}
\end{figure}

Below, in Section \ref{sec:numerics}, we consider an explicit example function constructed from tensor products of B-splines (see also Figure \ref{fig:Bspline2:nonper}). 


\begin{remark}
For the purposes of this short note we have formulated the results of this section in the Hilbert space setting. However, with the same techniques we are able to prove analogous results for the more general Besov spaces $S^{r}_{p,q}B(D)$ of dominating mixed smoothness as well. The definition of these spaces is based on generalizations of the norm in equation (\ref{eq:besov_type_norm}).

Within the context of Besov spaces the operator $T_{\cos}$ can be interpreted not just as a periodization, but also as a generalized change of variables. In this regard our work is closely related to \cite{bourdaud2000, ikeda2023boundedness, NguyenDiss} and \cite{Nguyen2017}.

Finally, in analogy to the classical function spaces on the torus based on the Fourier system it makes sense to define function spaces based on Chebyshev coefficients. A theory of such spaces in an isotropic setting is developed in the works \cite{berthold1992fast, glaeske1987discrete, ivanov2012decomposition, junghanns2021weighted, kyriazis2008decomposition, kyriazis2008jacobi, runst1980strong, triebel2010theory} among others. Our results give novel insight into the relationship between the classical Sobolev spaces on $D$ and these Chebyshev--type spaces in the $d$-dimensional setting.
\end{remark}


\section{Subsampling} \label{sec:subsampling}

To efficiently recover functions on $D$ we need to find suitable sample points $\bx^i \in D, i \in [n]$. This class of points will be used for the whole function space and thus needs to be calculated only once. To do so, we start by sampling $\mathcal{O}(m \log m)$ random points according to the probability measure $\varrho^{}_D$ from \eqref{eq:rho} and use a subsampling technique based on Weaver's $KS_2$-conjecture (proven in \cite{MSS15} and further developed in \cite{NOU15}, see \cite{nagelBA} for a complete account) to reduce the budget to $\mathcal{O}(m)$. We start with the following theorem from \cite{NSU22,NOU15}.

\smallskip

\begin{theorem} \label{thm:subsampling_existence}
    Let $c, A, B > 0$ and $\bu^1, \dots, \bu^M \in \C^m$ with $\|\bu^i\|_2^2 \leq c \frac{m}{M}$ for all $i=1, \ldots, M$ and
    $$
    A \|\bw\|_2^2 \leq \sum_{i=1}^M |\langle \bw, \bu^i\rangle|^2 \leq B \|\bw\|_2^2
    $$
    for all $\bw \in \C^m$. Then there is a $J \subset [M]$ of size $\# J \leq bm$ with
    $$
    A' \cdot \frac{m}{M} \|\bw\|_2^2 \leq \sum_{j\in J} |\langle \bw, \bu^j\rangle|^2 \leq B' \cdot \frac{m}{M} \|\bw\|_2^2
    $$
    for all $\bw \in \C^m$, where $b, A', B' > 0$ only depend on $c, A, B$ (and can be given explicitly).
\end{theorem}

\smallskip

In Theorem \ref{thm:subsampling_existence} the constants $A$ and $B$ are called the frame bounds and the system $\bu^1, \dots, \bu^M$ is called a frame. The strategy here is to use $\bu^i = [\eta_\bk(\bx^i)]_{\bk \in \Lambda}$ as the frame, where $\Lambda$ is a set of frequencies. The indices $J \subset [M]$ as chosen in the above theorem will give sample points $\{\bx^j\}_{j \in J}$ suitable for approximation.

Theorem \ref{thm:subsampling_existence} only shows the existence of such a good set of sampling points. It also leaves no control on the oversampling factor $b$ (e.g. if $c=1$, $A=1/2$ and $B=3/2$, Theorem \ref{thm:subsampling_existence} gives $b=3284$, $A'=11.65$ and $B'=4926$). Constructive subsampling methods with arbitrary oversampling factors $b > 1 + \varepsilon$ were studied in \cite{BSU22}.

%

\smallskip

\begin{theorem} \label{thm:subsampling_constructive_unweighted}
    Let $\bu^1, ..., \bu^M \in \C^m$, choose $b > 1+\frac{1}{m}$ and assume $M \geq bm$. There is a polynomial time algorithm to construct a $J \subset [M]$ with $\# J \leq \lceil bm\rceil$ and
    $$
    \frac{1}{M} \sum_{i=1}^M |\langle \bw, \bu^i\rangle|^2 \leq 89 \frac{(b+1)^2}{(b-1)^3} \cdot \frac{1}{m} \sum_{j\in J} |\langle \bw, \bu^j\rangle|^2
    $$
    for all $\bw \in \C^m$.
\end{theorem}

\smallskip

Theorem \ref{thm:subsampling_constructive_unweighted} enables us to choose an oversampling factor $b$ close to $1$, but it also has the disadvantage of only giving a lower bound. However, this suffices for our purposes. We will use Algorithms 1 and 3 from \cite{BSU22} for our numerical experiments.

\section{Numerical Experiments} \label{sec:numerics}

\begin{figure}[htb]
		\begin{tikzpicture}[baseline]
		\begin{axis}[
			font=\footnotesize,
			enlarge x limits=true,
			enlarge y limits=true,
			height=0.4\textwidth,
			grid=major,
			width=0.47\textwidth,
            xmode=log,
            ymode=log,
			xlabel={$n$},
			ylabel={$L_2(\varrho^{}_D)$-error},
			legend style={legend cell align=right, at={(1.0,1.33)}},
			legend columns = 1,
		]
		\addplot[blue,mark=square,mark size=2.2pt,mark options={solid}] coordinates {
 (105, 0.002493838827094228) (132, 0.005470570995542974) (152, 0.0019317270825603948) (182, 0.0016201647863787292) (191, 0.0007342145619100965) (229, 0.0023438488410772113) (248, 0.0009296696785573356) (272, 0.0012200986579042911) (296, 0.0005276569023398316) (328, 0.0005447346073840774) (346, 0.00044715935171410064) (384, 0.0003810166401789718) (397, 0.00027259307193623595) (442, 0.000398132189088048) (468, 0.0002785376850259007) (489, 0.0003077333841955345) (506, 0.00017534258017248605) (569, 0.0003068486395050551) (584, 0.00015002502888240968) (756, 0.00010623240888879234) (883, 8.168633869702274e-05) (1092, 6.0824462145599186e-05) (1257, 3.078573880245972e-05) (1449, 2.7126693197178757e-05)
};
\addlegendentry{$d=3$}
\addplot [forget plot,black,domain=300:1600, samples=100, dotted,ultra thick]{0.1*x^(-2.5)*(ln(x))^(2*2.5+0.5)};

		\addplot[magenta,mark=o,mark size=2.5pt,mark options={solid}] coordinates {
 (145, 0.0096515647393509) (179, 0.002291422805289902) (269, 0.004137409180526302) (301, 0.005194488603926028) (392, 0.005016345812050545) (459, 0.001154875647956181) (550, 0.0007383631879260684) (579, 0.0009249091924040521) (739, 0.0005462599898081653) (775, 0.0005686463685449755) (863, 0.0007809894318122451) (953, 0.00032460892824195064) (1083, 0.00035680476657611135) (1135, 0.0003047028433068013) (1289, 0.0002335401919999115) (1327, 0.0002369992790281783) (1477, 0.0002558698166443881) (1583, 0.0001727758946981366) (1666, 0.0001762984613367274) (1687, 0.00017142414058306273) (1976, 0.00013241734608910536) (2026, 0.00015004235847634835) (2652, 6.575055972225433e-05)
};
\addlegendentry{$d=4$}
\addplot [forget plot,black,domain=500:3000, samples=100, dotted,ultra thick]{0.004*x^(-2.5)*(ln(x))^(3*2.5+0.5)};

		\addplot[green,mark=triangle,mark size=3pt,mark options={solid,rotate=180}] coordinates {
 (119, 0.024839324149489188) (205, 0.006977876724297209) (381, 0.007531861524229203) (457, 0.005027097318372271) (722, 0.002598839421555864) (813, 0.0013240524562349028) (1117, 0.001461665327517249) (1274, 0.000730866993063157) (1642, 0.0004232867067081064) (2208, 0.0004727477298909333) (2539, 0.00037736125263974554) (2808, 0.00022604767081089457)
};
\addlegendentry{$d=5$}
\addplot [forget plot,black,domain=500:3000, samples=100, dotted,ultra thick]{0.00015*x^(-2.5)*(ln(x))^(4*2.5+0.5)};

\addlegendimage{black,dotted,ultra thick}
\addlegendentry{$\sim n^{-2.5} (\log n)^{2.5(d-1) + 0.5}$}
		\end{axis}
		\end{tikzpicture}
\caption{Sampling errors in the $L_2(\varrho^{}_D)$-norm using the Chebyshev basis for $f\in H^{5/2-\varepsilon}_{\textnormal{mix}}(D)$ from \eqref{eq:testfunction}.}
\label{fig:L2_error_cheb_subsampled}
\end{figure}
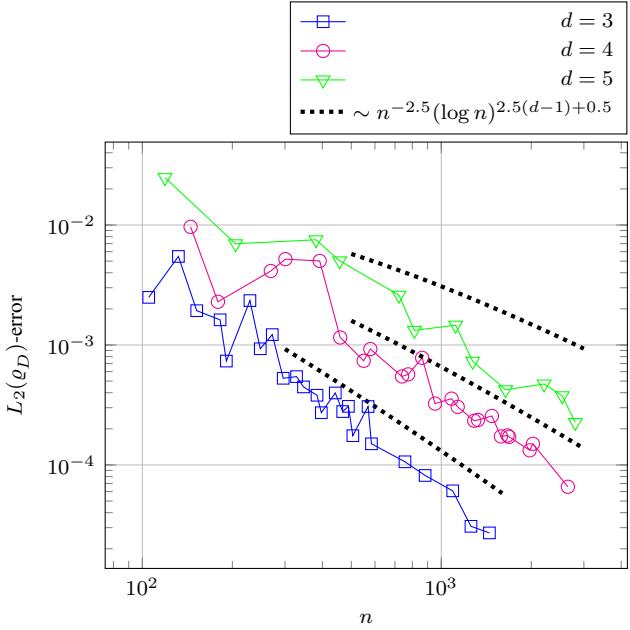

As an example, we follow \cite{Bartel22} and consider the function
\begin{align} \label{eq:testfunction}
\textstyle
    f(\bx) = \prod_{\ell=1}^{d} B_{2}(x_{\ell}) 
\end{align}
with
\begin{equation}
    B_{2}(x)
    = \begin{cases}
        -\frac{x^{2}}{4} -\frac{x}{2} + \frac{1}{2} &, \; -1 \leq x \leq 0 \\
        \frac{x^{2}}{8} -\frac{x}{2} + \frac{1}{2} &, \; 0 < x \leq 1
    \end{cases}. \notag
\end{equation}
$B_{2}$ is a cutout of a piecewise quadratic B-spline. Since spaces of dominating mixed smoothness have the cross norm property we observe that $f \in B^{5/2}_{2, \infty}(D)$ or, on the scale of Sobolev spaces, $f \in H^{s}_{\textnormal{mix}}(D)$ if and only if $s < 5/2$ (see \eqref{eq:chebcoeffs} below). This fact is reflected in the observed decay rates in Figure \ref{fig:L2_error_cheb_subsampled}. In terms of the $L_2(\varrho^{}_D)$-normalized Chebyshev polynomials of the first kind $T_{k}$ from \eqref{eq:cheb_poly} we have
\begin{align} \label{eq:chebcoeffs}
\begin{split}
    B_2(x) & = \frac{15}{32} T_0(x) - \frac{\pi-1}{2\pi\sqrt{2}} T_1(x) - \frac{1}{32\sqrt{2}} T_2(x) \\
    & - \frac{3}{2\pi\sqrt{2}} \sum_{k=3}^\infty \frac{ \sin(k\pi/2)}{k(k^2-4)} T_k(x).
\end{split}
\end{align}

\begin{figure}[htb]
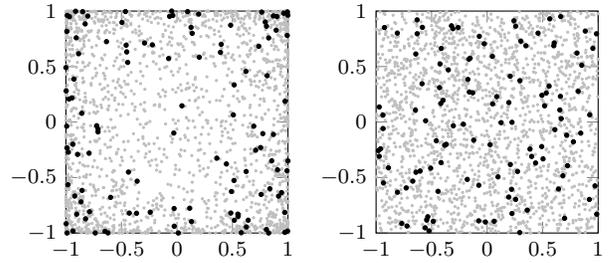

\centering

\caption{$M = 2000$ random Chebyshev [uniform] sample points (light gray, $A \approx 0.766$ [$0.375$], $B \approx 1.225$ [$0.628$]), subsampled to $n = 117$ [$116$] points (black, $A \approx 0.026$ [$0.01$], $B \approx 1.872$ [$1.068$]), using $m = 107$ frequencies (hyperbolic cross, oversampling factor $b = 1.1$)}
\label{fig:subsampling}
\end{figure}

We approximate functions over $D$ via the following procedure. Given a hyperbolic cross
$$
\Lambda = \Lambda_{d, R} = \{\bk \in \N_0^d: \prod_{\ell=1}^d \max\{1, k_\ell\} \leq R\}
$$
with $d, R \in \N$, set $m \coloneqq \#\Lambda$ (for estimates on the size of $\Lambda_{d, R}$ see \cite{CHERNOV201692} and \cite{KSU13}, but note that the hyperbolic cross is defined slightly different there). We now draw $M \coloneqq \lceil 4m \log m \rceil$ random points according to the probability measure $\varrho^{}_D$ (see Figure \ref{fig:subsampling} (left plot) and note the higher concentration of points at the boundary). Letting this set be denoted by $X = \{\bx^i: i \in [M]\}$, we define the vectors $\bu^i = \left[\eta_\bk(\bx^i)\right]_{\bk \in \Lambda}$ with $\eta_\bk(\bx^i)$ as in \eqref{eq:eta}. Using an oversampling factor of $b=1.1$ the algorithm from Theorem \ref{thm:subsampling_constructive_unweighted} gives a set $J = \{j_1, \dots, j_n\} \subset [M]$ of size $\#J = n \leq \lceil bm\rceil$, such that the $\mathcal{O}(m)$ points $\{\bx^j\}_{j \in J}$ approximate functions on $[-1, 1]^d$ just as well (asymptotically) as the $\mathcal{O}(m \log m)$ random points $X$.

Figure \ref{fig:subsampling} (left plot) shows the subsampling step for $d = 2$ and $R = 20$, so that $m = 107$ and $M = 2000$. The algorithm finds a subset of $n = 117 \leq \lceil 1.1 \cdot 107\rceil$ points that are suitable for function approximation. The quality of a point set in this case is measured by the smallest and largest singular values (i.e. the condition) of the matrix $(\# X)^{-1/2} \left[\eta_\bk(\bx)\right]_{\bx \in X, \bk \in \Lambda}$. For the $M = 2000$ random points the smallest and largest singular values respectively are $A \approx 0.766$ and $B \approx 1.225$. For the $n = 117$ subsampled points these come out as $A \approx 0.026$ and $B \approx 1.872$. 

Calculating coefficients $c_\bk, \bk \in \Lambda$ via the least squares system
$$
\left[\begin{matrix}
\eta_{\bk_1}(\bx^{j_1}) & \cdots & \eta_{\bk_m}(\bx^{j_1}) \\
\vdots & & \vdots \\
\eta_{\bk_1}(\bx^{j_n}) & \cdots & \eta_{\bk_m}(\bx^{j_n})
\end{matrix} \right] \left[\begin{matrix}
c_{\bk_1} \\
\vdots \\
c_{\bk_m}
\end{matrix}\right] \approx \left[\begin{matrix}
f(\bx^{j_1}) \\
\vdots \\
f(\bx^{j_n})
\end{matrix}\right]
$$
we approximate $f \approx \sum_{\bk \in \Lambda} c_\bk \eta_\bk$. Figure \ref{fig:L2_error_cheb_subsampled} shows the sampling error $\|f-\sum_{\bk \in \Lambda} c_\bk \eta_\bk\|_{L_2(\varrho^{}_D)}$ against the number of samples $n$ (after the subsampling step). Note that this procedure is inherently probabilistic, explaining the fluctuations in the plot. We compare this to the theoretical bound $n^{-2.5} (\log n)^{2.5(d-1)+0.5}$ from \cite{NSU22}, Section 7 (building on concentration inequalities for matrices from \cite{MoeUl20}) and improving over \cite{KaUlVo19}, see also \cite{TEMLYAKOV2021101545}. The authors in \cite{DKU22} give the final sharp, however  non-constructive, version based on Theorem \ref{thm:subsampling_existence}, which solved an outstanding open problem from \cite{DuTeUl2015}.

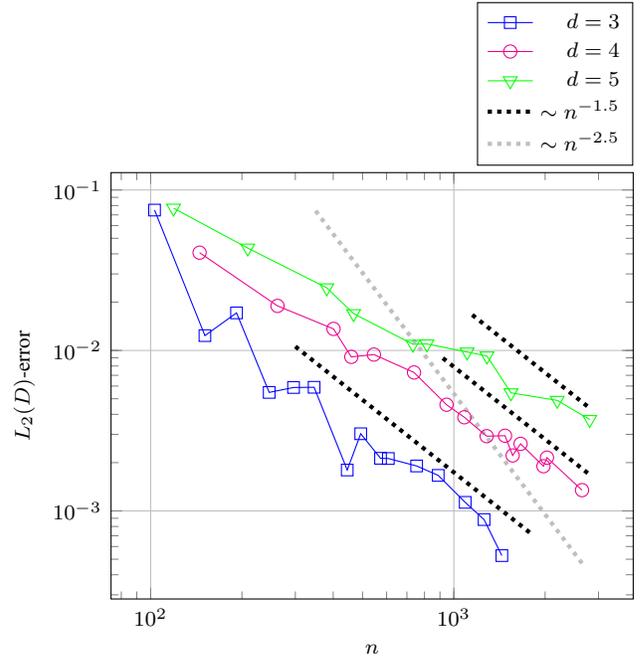
\begin{figure}[htb]
		\begin{tikzpicture}[baseline]
		\begin{axis}[
			font=\footnotesize,
			enlarge x limits=true,
			enlarge y limits=true,
			height=0.4\textwidth,
			grid=major,
			width=0.47\textwidth,
            xmode=log,
            ymode=log,
			xlabel={$n$},
			ylabel={$L_2(D)$-error},
			legend style={legend cell align=right, at={(1.0,1.4)}},
			legend columns = 1,
		]
		
\addplot [forget plot,lightgray,domain=350:2650, samples=100, dotted,ultra thick]{170000*x^(-2.5)};
		
		\addplot[blue,mark=square,mark size=2.2pt,mark options={solid}] coordinates {
 (103, 0.0749890561952042) (151, 0.012376903254038993) (192, 0.01714590622864491) (246, 0.005477199164133434) (296, 0.005877971540307987) (346, 0.005896470789471396) (445, 0.0017988374659264968) (493, 0.003026488591155279) (573, 0.002128442031775004) (608, 0.0021215397377644875) (754, 0.001906530091611939) (888, 0.0016650411945020136) (1090, 0.0011306938351054636) (1258, 0.0008831675903132739) (1440, 0.0005266321263162264)
};
\addlegendentry{$d=3$}
\addplot [forget plot,black,domain=300:1800, samples=100, dotted,ultra thick]{55 * x^(-1.5)};

		\addplot[magenta,mark=o,mark size=2.5pt,mark options={solid}] coordinates {
 (145, 0.04069585400777723) (262, 0.018962434703360214) (401, 0.013627197499881841) (459, 0.009119532595600327) (545, 0.009427207847427154) (739, 0.00729090117268837) (946, 0.004590237587796813) (1083, 0.0038346101838748) (1283, 0.002923623372703624) (1475, 0.0029423923754519085) (1563, 0.002216753668075339) (1660, 0.0026182880105047003) (1975, 0.0018994854644859482) (2028, 0.0021485276754043075) (2647, 0.0013473637753040563)
};
\addlegendentry{$d=4$}
\addplot [forget plot,black,domain=920:2800, samples=100, dotted,ultra thick]{250 * x^(-1.5)};

		\addplot[green,mark=triangle,mark size=3pt,mark options={solid,rotate=180}] coordinates {
(119, 0.07708949117869114) (209, 0.04346058248894587) (381, 0.02455675605462107) (467, 0.01695825944729583) (733, 0.010934003173981407) (815, 0.010979844908341212) (1106, 0.009774812018157492) (1286, 0.009236866895767814) (1543, 0.005445835804773233) (2193, 0.004885065657104202) (2806, 0.003711678105891124)
};
\addlegendentry{$d=5$}
\addplot [forget plot,black,domain=1150:2800, samples=100, dotted,ultra thick]{650 * x^(-1.5)};

\addlegendimage{black,dotted,ultra thick}
\addlegendentry{$\sim n^{-1.5}$}
\addlegendimage{lightgray,dotted,ultra thick}
\addlegendentry{$\sim n^{-2.5}$}
		\end{axis}
		\end{tikzpicture}
\caption{Sampling errors in the $L_2(D)$-norm using the half period cosine basis for $f\in H^{5/2-\varepsilon}_{\textnormal{mix}}(D)$ from \eqref{eq:testfunction}.}
\label{fig:L2_error_hpc_subsampled}
\end{figure}

This gives an improved rate compared to the approximation with the half period cosine basis $V_k(x) = \sqrt{2}^{-\delta_{0, k}} \cos\left(\pi k \frac{x+1}{2}\right)$ ($L_2(D)$-normalized) over $[-1, 1]$, where we have
$$
B_2(x) = \frac{23}{24\sqrt{2}} V_0(x) + \sum_{k=1}^\infty \left(\frac{6\sin(\pi k/2)}{\pi^3 k^3} - \frac{(-1)^{k}}{\pi^2 k^2}\right) V_k(x).
$$
Compared to the above, we only need to take uniform samples over $[-1, 1]^d$ and change the basis functions to $\psi_\bk(\bx^i) = \prod_{\ell=1}^d V_{k_\ell}(x^i_\ell)$. Figure \ref{fig:subsampling} (right plot) demonstrates the subsampling step and Figure \ref{fig:L2_error_hpc_subsampled} shows the sampling error (with respect to $L_2(D)$). While other, smaller order factors (e.g. polylogarithmic) are not known, we still observe a main decay rate of $3/2$ (compared to $5/2$ in the Chebyshev setting).

Note that the half period cosine is the natural basis of $H^1_{\textnormal{mix}}(D)$ and the linear decay is guaranteed. In \cite{Bartel22}, the basis of $H^2_{\textnormal{mix}}(D)$ was considered with guaranteed quadratic decay. In the numerical experiment there the rate $5/2$ was obtained similar to the Chebyshev setting.

We also want to emphasize that the error in Figure \ref{fig:L2_error_cheb_subsampled} is measured in the stronger $L_2(\varrho^{}_D)$-norm compared to the $L_2(D)$-norm used in Figure \ref{fig:L2_error_hpc_subsampled}. This is due to
$$
\int_D |g(\bx)|^2 \mathrm{d}\varrho^{}_D(\bx) \geq \pi^{-d} \int_D |g(\bx)|^2 \mathrm{d}\bx
$$
for functions $g$, for which the integrals are defined (as is the case for continuous functions on $D$), while no comparable bound can be made in the other direction.

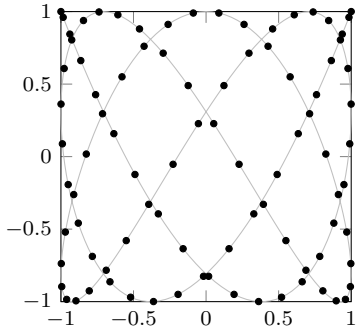
\begin{figure}[htb]
\centering
		\begin{tikzpicture}
		\begin{axis}[
		plot box ratio = 1 1,
		xmin=-1,xmax=1,ymin=-1,ymax=1,
		font=\footnotesize,
		height=0.3\textwidth,
		width=0.3\textwidth,
		]
		\addplot[lightgray] coordinates {
			(1.0, 1.0) (0.999, 0.9975) (0.9961, 0.99) (0.9913, 0.9777) (0.9845, 0.9604) (0.9758, 0.9383) (0.9652, 0.9116) (0.9527, 0.8804) (0.9383, 0.8447) (0.9222, 0.8048) (0.9042, 0.761) (0.8845, 0.7133) (0.8631, 0.6621) (0.8399, 0.6076) (0.8152, 0.5501) (0.7889, 0.4898) (0.761, 0.4271) (0.7316, 0.3622) (0.7008, 0.2956) (0.6687, 0.2274) (0.6352, 0.1582) (0.6006, 0.08813) (0.5647, 0.01765) (0.5278, -0.05292) (0.4898, -0.1232) (0.4508, -0.1929) (0.411, -0.2617) (0.3704, -0.3291) (0.3291, -0.3949) (0.2871, -0.4587) (0.2446, -0.5202) (0.2016, -0.5792) (0.1582, -0.6352) (0.1145, -0.6881) (0.07054, -0.7376) (0.02647, -0.7834) (-0.01765, -0.8253) (-0.06173, -0.8631) (-0.1057, -0.8965) (-0.1495, -0.9256) (-0.1929, -0.95) (-0.236, -0.9696) (-0.2787, -0.9845) (-0.3207, -0.9944) (-0.3622, -0.9994) (-0.403, -0.9994) (-0.4429, -0.9944) (-0.4821, -0.9845) (-0.5202, -0.9696) (-0.5574, -0.95) (-0.5935, -0.9256) (-0.6284, -0.8965) (-0.6621, -0.8631) (-0.6945, -0.8253) (-0.7256, -0.7834) (-0.7552, -0.7376) (-0.7834, -0.6881) (-0.8101, -0.6352) (-0.8351, -0.5792) (-0.8586, -0.5202) (-0.8804, -0.4587) (-0.9004, -0.3949) (-0.9187, -0.3291) (-0.9353, -0.2617) (-0.95, -0.1929) (-0.9628, -0.1232) (-0.9738, -0.05292) (-0.9829, 0.01765) (-0.99, 0.08813) (-0.9953, 0.1582) (-0.9986, 0.2274) (-1.0, 0.2956) (-0.9994, 0.3622) (-0.9968, 0.4271) (-0.9924, 0.4898) (-0.986, 0.5501) (-0.9777, 0.6076) (-0.9674, 0.6621) (-0.9553, 0.7133) (-0.9414, 0.761) (-0.9256, 0.8048) (-0.908, 0.8447) (-0.8886, 0.8804) (-0.8675, 0.9116) (-0.8447, 0.9383) (-0.8203, 0.9604) (-0.7942, 0.9777) (-0.7667, 0.99) (-0.7376, 0.9975) (-0.7071, 1.0) (-0.6752, 0.9975) (-0.642, 0.99) (-0.6076, 0.9777) (-0.572, 0.9604) (-0.5352, 0.9383) (-0.4975, 0.9116) (-0.4587, 0.8804) (-0.4191, 0.8447) (-0.3786, 0.8048) (-0.3374, 0.761) (-0.2956, 0.7133) (-0.2531, 0.6621) (-0.2102, 0.6076) (-0.1669, 0.5501) (-0.1232, 0.4898) (-0.07934, 0.4271) (-0.03529, 0.3622) (0.008825, 0.2956) (0.05292, 0.2274) (0.09692, 0.1582) (0.1407, 0.08813) (0.1843, 0.01765) (0.2274, -0.05292) (0.2702, -0.1232) (0.3124, -0.1929) (0.354, -0.2617) (0.3949, -0.3291) (0.435, -0.3949) (0.4743, -0.4587) (0.5127, -0.5202) (0.5501, -0.5792) (0.5864, -0.6352) (0.6215, -0.6881) (0.6555, -0.7376) (0.6881, -0.7834) (0.7195, -0.8253) (0.7494, -0.8631) (0.7779, -0.8965) (0.8048, -0.9256) (0.8302, -0.95) (0.854, -0.9696) (0.8761, -0.9845) (0.8965, -0.9944) (0.9152, -0.9994) (0.9321, -0.9994) (0.9472, -0.9944) (0.9604, -0.9845) (0.9717, -0.9696) (0.9812, -0.95) (0.9888, -0.9256) (0.9944, -0.8965) (0.9981, -0.8631) (0.9998, -0.8253) (0.9996, -0.7834) (0.9975, -0.7376) (0.9934, -0.6881) (0.9874, -0.6352) (0.9795, -0.5792) (0.9696, -0.5202) (0.9579, -0.4587) (0.9443, -0.3949) (0.9289, -0.3291) (0.9116, -0.2617) (0.8926, -0.1929) (0.8718, -0.1232) (0.8494, -0.05292) (0.8253, 0.01765) (0.7996, 0.08813) (0.7723, 0.1582) (0.7435, 0.2274) (0.7133, 0.2956) (0.6817, 0.3622) (0.6488, 0.4271) (0.6146, 0.4898) (0.5792, 0.5501) (0.5427, 0.6076) (0.5051, 0.6621) (0.4665, 0.7133) (0.4271, 0.761) (0.3868, 0.8048) (0.3457, 0.8447) (0.304, 0.8804) (0.2617, 0.9116) (0.2188, 0.9383) (0.1756, 0.9604) (0.132, 0.9777) (0.08813, 0.99) (0.04411, 0.9975) (1.194e-15, 1.0) (-0.04411, 0.9975) (-0.08813, 0.99) (-0.132, 0.9777) (-0.1756, 0.9604) (-0.2188, 0.9383) (-0.2617, 0.9116) (-0.304, 0.8804) (-0.3457, 0.8447) (-0.3868, 0.8048) (-0.4271, 0.761) (-0.4665, 0.7133) (-0.5051, 0.6621) (-0.5427, 0.6076) (-0.5792, 0.5501) (-0.6146, 0.4898) (-0.6488, 0.4271) (-0.6817, 0.3622) (-0.7133, 0.2956) (-0.7435, 0.2274) (-0.7723, 0.1582) (-0.7996, 0.08813) (-0.8253, 0.01765) (-0.8494, -0.05292) (-0.8718, -0.1232) (-0.8926, -0.1929) (-0.9116, -0.2617) (-0.9289, -0.3291) (-0.9443, -0.3949) (-0.9579, -0.4587) (-0.9696, -0.5202) (-0.9795, -0.5792) (-0.9874, -0.6352) (-0.9934, -0.6881) (-0.9975, -0.7376) (-0.9996, -0.7834) (-0.9998, -0.8253) (-0.9981, -0.8631) (-0.9944, -0.8965) (-0.9888, -0.9256) (-0.9812, -0.95) (-0.9717, -0.9696) (-0.9604, -0.9845) (-0.9472, -0.9944) (-0.9321, -0.9994) (-0.9152, -0.9994) (-0.8965, -0.9944) (-0.8761, -0.9845) (-0.854, -0.9696) (-0.8302, -0.95) (-0.8048, -0.9256) (-0.7779, -0.8965) (-0.7494, -0.8631) (-0.7195, -0.8253) (-0.6881, -0.7834) (-0.6555, -0.7376) (-0.6215, -0.6881) (-0.5864, -0.6352) (-0.5501, -0.5792) (-0.5127, -0.5202) (-0.4743, -0.4587) (-0.435, -0.3949) (-0.3949, -0.3291) (-0.354, -0.2617) (-0.3124, -0.1929) (-0.2702, -0.1232) (-0.2274, -0.05292) (-0.1843, 0.01765) (-0.1407, 0.08813) (-0.09692, 0.1582) (-0.05292, 0.2274) (-0.008825, 0.2956) (0.03529, 0.3622) (0.07934, 0.4271) (0.1232, 0.4898) (0.1669, 0.5501) (0.2102, 0.6076) (0.2531, 0.6621) (0.2956, 0.7133) (0.3374, 0.761) (0.3786, 0.8048) (0.4191, 0.8447) (0.4587, 0.8804) (0.4975, 0.9116) (0.5352, 0.9383) (0.572, 0.9604) (0.6076, 0.9777) (0.642, 0.99) (0.6752, 0.9975) (0.7071, 1.0) (0.7376, 0.9975) (0.7667, 0.99) (0.7942, 0.9777) (0.8203, 0.9604) (0.8447, 0.9383) (0.8675, 0.9116) (0.8886, 0.8804) (0.908, 0.8447) (0.9256, 0.8048) (0.9414, 0.761) (0.9553, 0.7133) (0.9674, 0.6621) (0.9777, 0.6076) (0.986, 0.5501) (0.9924, 0.4898) (0.9968, 0.4271) (0.9994, 0.3622) (1.0, 0.2956) (0.9986, 0.2274) (0.9953, 0.1582) (0.99, 0.08813) (0.9829, 0.01765) (0.9738, -0.05292) (0.9628, -0.1232) (0.95, -0.1929) (0.9353, -0.2617) (0.9187, -0.3291) (0.9004, -0.3949) (0.8804, -0.4587) (0.8586, -0.5202) (0.8351, -0.5792) (0.8101, -0.6352) (0.7834, -0.6881) (0.7552, -0.7376) (0.7256, -0.7834) (0.6945, -0.8253) (0.6621, -0.8631) (0.6284, -0.8965) (0.5935, -0.9256) (0.5574, -0.95) (0.5202, -0.9696) (0.4821, -0.9845) (0.4429, -0.9944) (0.403, -0.9994) (0.3622, -0.9994) (0.3207, -0.9944) (0.2787, -0.9845) (0.236, -0.9696) (0.1929, -0.95) (0.1495, -0.9256) (0.1057, -0.8965) (0.06173, -0.8631) (0.01765, -0.8253) (-0.02647, -0.7834) (-0.07054, -0.7376) (-0.1145, -0.6881) (-0.1582, -0.6352) (-0.2016, -0.5792) (-0.2446, -0.5202) (-0.2871, -0.4587) (-0.3291, -0.3949) (-0.3704, -0.3291) (-0.411, -0.2617) (-0.4508, -0.1929) (-0.4898, -0.1232) (-0.5278, -0.05292) (-0.5647, 0.01765) (-0.6006, 0.08813) (-0.6352, 0.1582) (-0.6687, 0.2274) (-0.7008, 0.2956) (-0.7316, 0.3622) (-0.761, 0.4271) (-0.7889, 0.4898) (-0.8152, 0.5501) (-0.8399, 0.6076) (-0.8631, 0.6621) (-0.8845, 0.7133) (-0.9042, 0.761) (-0.9222, 0.8048) (-0.9383, 0.8447) (-0.9527, 0.8804) (-0.9652, 0.9116) (-0.9758, 0.9383) (-0.9845, 0.9604) (-0.9913, 0.9777) (-0.9961, 0.99) (-0.999, 0.9975) (-1.0, 1.0)
		};
		\addplot[only marks,black,mark=*,mark size=1.2pt,mark options={solid}] coordinates {
		    (1.0, 1.0) (0.9845, 0.9604) (0.9383, 0.8447) (0.8631, 0.6621) (0.761, 0.4271) (0.6352, 0.1582) (0.4898, -0.1232) (0.3291, -0.3949) (0.1582, -0.6352) (-0.01765, -0.8253) (-0.1929, -0.95) (-0.3622, -0.9994) (-0.5202, -0.9696) (-0.6621, -0.8631) (-0.7834, -0.6881) (-0.8804, -0.4587) (-0.95, -0.1929) (-0.99, 0.08813) (-0.9994, 0.3622) (-0.9777, 0.6076) (-0.9256, 0.8048) (-0.8447, 0.9383) (-0.7376, 0.9975) (-0.6076, 0.9777) (-0.4587, 0.8804) (-0.2956, 0.7133) (-0.1232, 0.4898) (0.05292, 0.2274) (0.2274, -0.05292) (0.3949, -0.3291) (0.5501, -0.5792) (0.6881, -0.7834) (0.8048, -0.9256) (0.8965, -0.9944) (0.9604, -0.9845) (0.9944, -0.8965) (0.9975, -0.7376) (0.9696, -0.5202) (0.9116, -0.2617) (0.8253, 0.01765) (0.7133, 0.2956) (0.5792, 0.5501) (0.4271, 0.761) (0.2617, 0.9116) (0.08813, 0.99) (-0.08813, 0.99) (-0.2617, 0.9116) (-0.4271, 0.761) (-0.5792, 0.5501) (-0.7133, 0.2956) (-0.8253, 0.01765) (-0.9116, -0.2617) (-0.9696, -0.5202) (-0.9975, -0.7376) (-0.9944, -0.8965) (-0.9604, -0.9845) (-0.8965, -0.9944) (-0.8048, -0.9256) (-0.6881, -0.7834) (-0.5501, -0.5792) (-0.3949, -0.3291) (-0.2274, -0.05292) (-0.05292, 0.2274) (0.1232, 0.4898) (0.2956, 0.7133) (0.4587, 0.8804) (0.6076, 0.9777) (0.7376, 0.9975) (0.8447, 0.9383) (0.9256, 0.8048) (0.9777, 0.6076) (0.9994, 0.3622) (0.99, 0.08813) (0.95, -0.1929) (0.8804, -0.4587) (0.7834, -0.6881) (0.6621, -0.8631) (0.5202, -0.9696) (0.3622, -0.9994) (0.1929, -0.95) (0.01765, -0.8253) (-0.1582, -0.6352) (-0.3291, -0.3949) (-0.4898, -0.1232) (-0.6352, 0.1582) (-0.761, 0.4271) (-0.8631, 0.6621) (-0.9383, 0.8447) (-0.9845, 0.9604) (-1.0, 1.0)
		};
		\end{axis}
		\end{tikzpicture}
\caption{Rank-1 Chebyshev lattice based on a Fibonacci lattice}
\label{fig:r1c_lattice}
\end{figure}

Note that our procedure of finding the appropriate sampling nodes is probabilistic in nature. However, once the node set is fixed it works for the whole class of functions (unit ball in a certain Sobolev space) unlike the Monte-Carlo methods like in \cite{DOLBEAULT2022101602}, where the node set might change for each function to approximate. One may ask for deterministic constructions of point sets to approximate non-periodic functions over $D$. In \cite{PV15} and \cite{volkmerdiss} so called rank-1 Chebyshev lattices (see Figure \ref{fig:r1c_lattice}) were considered. In the same way as above, one uses function evaluations at those lattice points to get the coefficients for tensored Chebyshev basis functions $\eta_\bk$ at frequencies given by a hyperbolic cross. Example 3.24 from \cite{volkmerdiss} also uses a cutout of a quadratic B-spline. While the theory there only gives a decay rate of $n^{-(1-\varepsilon)}$ for all $\varepsilon > 0$ (note the difference in notation there), the plots suggest an actual rate of $n^{-(1.25-\varepsilon)}$, which can be verified by our analysis. These constructive methods therefore only give half the rate compared to our method. Still, it should be noted that on lattices one can use FFT-techniques to accelerate the computation which are not possible for our unstructured (random) sample points. Also see \cite{barthelmann2000high} for further methods using sparse grids.

\section{Conclusion}

We give analytic and numerical evidence that for the recovery of non-periodic functions belonging to classical spaces over $[-1, 1]^d$ sampling points following a Chebyshev distribution together with Chebyshev polynomials (see Figure \ref{fig:L2_error_cheb_subsampled}) achieve state of the art error decay. Section \ref{sec:spaces} gives the theoretical framework for this where the periodization operator $T_{\cos}$ plays a crucial role in the analysis. We conclude that transforming non-periodic functions to periodic ones using the cosine composition $T_{\cos}$ helps in finding efficient approximation algorithms.

\section*{Acknowledgement}

The authors would like to thank Winfried Sickel and Erich Novak for suggesting further references. Felix Bartel would like to thank the Deutscher Akademischer Austauschdienst (DAAD) for funding his research scholarship. Nicolas Nagel is funded by the ESF ReSIDA-H2 project. Kai Lüttgen is supported by a \textit{Saxon Scholarship} granted by the Free State of Saxony.

\bibliographystyle{amsplain}
\bibliography{main.bib}

\end{document}